\documentclass[12pt]{article}

\usepackage{amssymb,amsmath,amsfonts}
\usepackage[T2A]{fontenc}
\usepackage[utf8]{inputenc}
\usepackage[english,russian]{babel}
\usepackage{graphicx}
\usepackage{epstopdf}

\let\chapter\section
\begin{document}

\begin{otherlanguage}{english}

{\large
\begin{center}
\textbf{Reducing smooth functions to normal forms near critical points}\\
A.S.~Orevkova
\end{center}
}

{\footnotesize
\noindent
s15b3\_orevkova@179.ru; Moscow State Univercity, Moscow Center of Fundamental and Applied Mathematics
}

\begin{abstract}
The paper is devoted to ``uniform'' reduction of smooth functions on 2-manifolds to canonical form near critical points by some coordinate changes in some neighbourhoods of these points. For singularity types $E_6,E_8$ and $A_n$, we explicitly construct such coordinate changes and estimate from below (in terms of $C^r$-norm of the function) the radius of a required neighbourhood.
\end{abstract}
\noindent
\textit{Keywords:} right equivalence of smooth functions, ADE-singularities, normal form of singularities, uniform reducing to normal form.

\noindent
\textit{MSC:} 58K05, 37C15.

\vspace{1em}
\hrule
\vspace{2em}

\section{Introduction}

\textbf{Definition.} {\it 
A smooth function $f = f(u_1, u_2) $ has a singularity type $E_k$ {\rm($k=6,7,8$)} at its critical point $P\in\mathbb {R}^2$ if
\begin{itemize}
\item[\rm(i)] the first and second differentials $df(P) = 0$, $d^2f(P)=0$, and the third differential $d^3f(P) \ne 0$ and is a perfect cube;
\item[\rm(ii)] one of the following coefficients of the Taylor series of $f$ at $P$ does not vanish: $f^{(4)}_{y^4}(P)$, $f^{(4)}_{xy^3}(P)$ and $f^{(5)}_{y^5}(P)$, where $(u_1,u_2)\to(x,y)$ is a linear coordinate change such that $d^3f(P)=(dx)^3$. More specifically: 
the singularity type is $E_6$ if $f^{(4)}_{y^4}(P)\ne0$ (equivalently, there exists a tangent vector $v\in \operatorname{Ker} d^3f(P)$ at $P$ such that $v^4f \ne 0$, where $v^4f$ denotes the fourth derivative of $f$ along the vector $ v $); 
the singularity type is $E_7$ if $f^{(4)}_{y^4}(P)=0$ and $f^{(4)}_{xy^3}(P)\ne 0$;
the singularity type is $E_8$ if $f^{(4)}_{y^4}(P)=f^{(4)}_{xy^3}(P)= 0$ and $f^{(5)}_{y^5}(P)\ne0$.
\end{itemize}}

From definition of $E_k$ we have $f^{(a+b)}_{x^ay^b}(P)=0$, 
$0<\frac{a}{3}+\frac{2b}{k+2}<1$, $f'''_{x^3}(P) = 6$ ($k=6,7,8$) and $f_{y^{1+k/2}}^{(1+k/2)}(P) \ne 0$ ($k=6,8$).
We will assume that $P={\bf 0}=(0,0)$ in the coordinates $x, y$.

\medskip
\textbf{Assumption 1.} For singularities $E_k$ ($k=6,8$), assume that
$f_{y^{1+k/2}}^{(1+k/2)}({\bf 0}) = \pm(1+\frac k2)!$.

\medskip
\textbf{Theorem 1} (Reducing $E_k$ to normal form {[1]})\textbf{.} 
{\it 
Let a function $ f (u_1, u_2) $ have a singularity $ E_k $ {\rm($k=6,7,8$)} at a critical point $ P $. Then, in some neighborhood of $P$, there is a local coordinate system $\tilde x, \tilde y$ in which the point $ P $ is the origin, and the function has the normal form $ f = f(P) + \tilde x^3 \pm \tilde y^{1+k/2}$ for $k=6,8$, $f = f(P) + \tilde x^3 +\tilde x \tilde y^3$ for $k=7$.}

\medskip 
In [1], the existence of a coordinate change was proved using the Tougeron theorem [13]. In view of this, obtaining a formula for the corresponding coordinate change requires solving the Cauchy problem for a system of ODE's. We construct our coordinate change explicitly, without using the Tougeron theorem.

\medskip
\textbf{Lemma 1} ([1])\textbf{.} {\it 
In the case of a singularity $E_k$ {\rm($k=6,8$)}, under the hypotheses of Assumption 1,
there exist coordinate changes 
${(x,y)} \to (x_1=x+d_1y^2,y) \to 
(x_1,y_1=y+d_2x_1) \to (x_2=x_1+d_3y_1^{k/2-1}, y_1) \to ({x_2, y_2}=y_1+d_4x_2^2)$ where
$d_j\in\mathbb R$ and $d_1=0$ for $k=6$, such that $f^{(a+b)}_{x_1^ay_2^b}({\bf 0})=0$ for all $a,b\in\mathbb{Z}_+$ with $a<3$ and $b<1+\frac k2$.}

\medskip
\textbf{Theorem 2} (Estimating the radius of a neighborhood for the coordinate change)\textbf{.} 
{\it In the case of singularities $E_k$ {\rm($k=6,8$)}, under the hypotheses of Assumption 1, let 
$(x,y)\to(x_2,y_2)$ be the coordinate change from Lemma 1. Suppose that, in a 
neighborhood $U_0 =\{(x_2,y_2)\mid \max(|x_2|,|y_2|) < R_0\}$ of $\bf 0$, the following estimates hold:
$C_{\alpha\beta}=\sup_{U_0}\big|f_{x_2^{\alpha}y_2^{\beta}}^{(\alpha+\beta)}(x_2,y_2)\big|\le M$ for $(\alpha,\beta)\in \{(0,5),(1,4),(3,1),(3,2),(3,3),(4,0),(4,1),(4,2),(4,3)\}$ if $k=6$, for $(\alpha,\beta)\in \{(0,6),(1,5),(3,1),(3,2),
(3,3),(3,4),(4,0),(4,1),(4,2),(4,3),(4,4)\}$ if $k=8$, where $R_0>0$, $M\ge0$. 
Then, in the neighbourhood $ U = \{(x_2, y_2) \mid \max (| x_2 |, | y_2 |) <R \} $, with $R =\min\{R_0,\frac{2}{M+2}\}$, there is a coordinate change of the form $\phi:(x_2,y_2)\to(\tilde x = x_2\sqrt[3]{h(x_2,y_2)}, \tilde y=y_2(g(x_2,y_2))^{\frac2{2+k}})$ 
that reduces $ f $ to the normal form $ f = f(P) + \tilde x^3 \pm \tilde y^{1+k/2}$ of $ E_k $. In detail:
\begin{itemize}
\item[\rm(a)] the functions $ h (\boldsymbol {x}) $ and $ g (\boldsymbol {x}) $ are positive in $ U $, thus the change $ \phi |_{U} $ is well-defined and is {$C^\infty$-smooth}; 
\item[\rm(b)] $\| \phi '(\boldsymbol {x}) - I \| <C <1 $ for all $\boldsymbol {x} \in U$, where $ C = \frac25 $, thus 
$ \phi |_{U} $ is $C^1$-close to the identity; 
\item[\rm(c)] the coordinate change $ \phi |_{U} $ is {injective} and {regular},
i.e., it is an embedding and $ \det | \phi '(\boldsymbol {x}) | \ne 0 $ for all $ \boldsymbol {x} \in U $, moreover $\phi(U)$

contains the open disk of radius $ (1 - C) R $ centred at ${\bf 0}$.
\end{itemize}}

\medskip
Our coordinate change $ \phi |_{U} $ from Theorem 2 provides a ``uniform'' reduction of the function $f$ at a singular point of type $ E_k $, $k=6,8$, to the canonical form $f=f(P)+ \tilde x ^ 3 \pm \tilde y ^ {1+k/2} $ in the sense that the neighbourhood radius and the coordinate change we constructed in this neighbourhood (as well as all partial derivatives of the coordinate change) continuously depend on the function $ f $ and its partial derivatives. 
A uniform reduction of smooth functions near critical points to a canonical form was known earlier for several singularity types [2, 10--12].

The uniform Morse lemma [3] was applied in [4--7] for studying topology of the spaces of Morse functions on surfaces and decomposition of these spaces into classes of topological equivalence. Our results have similar applications for stydying topology of the spaces of smooth functions [8] and gradient-like flows [9] with prescribed $ADE$-singularities.

The author is grateful to Elena Kudryavtseva for stating the problem and useful discussions.
The author is a Fellow of the Theoretical Physics and Mathematics Advancement Foundation ``BASIS''.

\section
{Key lemmas}

\medskip
\textbf{Lemma 2.} {\it Let $f^{(k+l)}_{x^ky^l}({\bf 0})=0$ for all $k,l\in\mathbb Z_+$ with $k<m$ and $l<n$, where $m,n\in\mathbb N$. Then the function $f|_U$ has the form $f(x,y)=x^mh(x,y) + y^ng(x,y)$ for some functions $h,g\in C^{\infty}(U)$, where $U$ is a neighborhood of $\bf 0$ in $\mathbb R^2$.}

\medskip
{\it Proof.} We represent $f$ as a function of $y$ with a parameter $x$. We write down the Taylor formula at the point $y = 0$ with the remainder term in the integral form:
$$
f(x,y) = S_n(x,y) + R_n(x,y) = x^mh(x,y) + y^ng(x,y),  
$$
$$
S_n(x,y)= f(x,0) + f_y'(x,0)y + f_{y^2}''(x,0)\frac{y^2}{2!} 
+ \dots + f_{y^{n-1}}^{(n-1)}(x,0)\frac{y^{n-1}}{(n-1)!}
=x^mh(x,y),
$$
$$
R_{n}(x,y) = \frac{1}{(n-1)!}\int_{0}^y f_{y^{n}}^{(n)}(x,t)(y-t)^{n-1} \,dt=
\frac{y^{n}}{(n-1)!}\int_{0}^1 f_{y^{n}}^{(n)}(x,sy)(1-s)^{n-1}\,ds = y^ng(x,y).
$$
Here each summand of $S_n(x,y)$ is the product of $y^i$ and a smooth function of $x$, which can be written as $x^mh_i(x)$, for some smooth functions $h_i(x)$. The lemma is proved. \hfill $\square$

\medskip
\textbf{Lemma 3.} {\it
Let $ \phi: U \to \mathbb R ^ n $ be a smooth mapping, where $ U $ is a convex open subset of $ \mathbb R ^ n $. Let the differential of $\phi$ have the form $ \phi'(\boldsymbol{x}) = I + A (\boldsymbol {x}) $, where $ I $ is the unit matrix of dimension $ n $, $ \| A \| < c $, $ 0 <c <1 $. Then $ \phi $ is injective and $ \det | \phi '(\boldsymbol{x}) | \ne0$ at every point $\boldsymbol{x} \in U $, i.e., $ \phi $ is a diffeomorphism to its image $ \phi (U) $. Moreover, $\langle \phi(\boldsymbol{x})-\phi(\boldsymbol{y}), \boldsymbol{x}-\boldsymbol{y} \rangle \ge (1-c) |\boldsymbol{x}-\boldsymbol{y}|^2$ for any pair of points $\boldsymbol{x},\boldsymbol{y}\in U$. }

\medskip

{\it Proof.}
Take any two points $ \boldsymbol {a}, \boldsymbol {b} \in U$, $\boldsymbol {a} \ne \boldsymbol {b} $, and consider the mapping $ \phi $ on the segment $ [\boldsymbol {a}, \boldsymbol {b}] $. The velocity vector
$$ 
v_t = \tfrac{d}{d t} \bigl (\phi(\boldsymbol {a} + t (\boldsymbol {b}-\boldsymbol {a})) \bigr). 
$$
Let us look at the projection of the velocity vector $ v_t $ onto $ \boldsymbol {b}-\boldsymbol {a} $:
\begin{equation*}
\begin{split}
\langle v_t,\boldsymbol {b}-\boldsymbol {a}\rangle &= \langle\tfrac{d}{d t}\left[\phi(\boldsymbol {a}+t(\boldsymbol {b}-\boldsymbol {a})) - \phi(\boldsymbol {a})\right],\boldsymbol {b}-\boldsymbol {a}\rangle = \langle \tfrac{d}{d t}\phi(\boldsymbol {a}+t(\boldsymbol {b}-\boldsymbol {a}))\rangle,\boldsymbol {b}-\boldsymbol {a}\rangle = \\
&= \langle \phi'(\boldsymbol {a}+t(\boldsymbol {b}-\boldsymbol {a}))(\boldsymbol {b}-\boldsymbol {a}),\boldsymbol {b}-\boldsymbol {a}\rangle = \langle (I+A(\boldsymbol {a}+t(\boldsymbol {b}-\boldsymbol {a})))\cdot(\boldsymbol {b}-\boldsymbol {a}),\boldsymbol {b}-\boldsymbol {a}\rangle = \\
&=\langle \boldsymbol {b}-\boldsymbol {a}, \boldsymbol {b}-\boldsymbol {a}\rangle + \langle A(\boldsymbol {a}+t(\boldsymbol {b}-\boldsymbol {a}))\cdot(\boldsymbol {b}-\boldsymbol {a}),(\boldsymbol {b}-\boldsymbol {a})\rangle.
\end{split}
\end{equation*}

Let us find an upper bound for the absolute value of the second term. Note that by the Cauchy-Schwarz inequality and by the definition of the matrix norm:
\begin{equation*}
\big|\langle A(\boldsymbol {a}+t(\boldsymbol {b}-\boldsymbol {a}))\cdot(\boldsymbol {b}-\boldsymbol {a}),(\boldsymbol {b}-\boldsymbol {a})\rangle\big| \le \|A(\boldsymbol {a}+t(\boldsymbol {b}-\boldsymbol {a}))\cdot(\boldsymbol {b}-\boldsymbol {a})\|\cdot\|\boldsymbol {b}-\boldsymbol {a}\| \le \|A\|\|\boldsymbol {b}-\boldsymbol {a}\|^2.
\end{equation*}

Let us go back to the estimation of $\langle v_t,\boldsymbol {b}-\boldsymbol {a}\rangle$:
\begin{equation*}
\begin{split}
\langle v_t,\boldsymbol {b}-\boldsymbol {a}\rangle &= \langle \boldsymbol {b}-\boldsymbol {a}, \boldsymbol {b}-\boldsymbol {a}\rangle + \langle A(\boldsymbol {a}+t(\boldsymbol {b}-\boldsymbol {a}))\cdot(\boldsymbol {b}-\boldsymbol {a}),\boldsymbol {b}-\boldsymbol {a}\rangle\geq \\ 
&\geq \|\boldsymbol {b}-\boldsymbol {a}\|^2 - |\langle A(\boldsymbol {a}+t(\boldsymbol {b}-\boldsymbol {a}))\cdot(\boldsymbol {b}-\boldsymbol {a}),\boldsymbol {b}-\boldsymbol {a}\rangle| \ge\\
&\geq \|\boldsymbol {b}-\boldsymbol {a}\|^2 - \|A\|\|\boldsymbol {b}-\boldsymbol {a}\|^2>\|\boldsymbol {b}-\boldsymbol {a}\|^2\cdot(1-c)>0,
\end{split}
\end{equation*}
from which $\det |\phi'(\boldsymbol {a})|\ne 0$.

By hypothesis $\|A\|<c<1$.
Let us look at the dot product $\langle \phi(\boldsymbol {b})-\phi(\boldsymbol {a}),\boldsymbol {b}-\boldsymbol {a}\rangle$. As $\langle \boldsymbol{x}(t), \boldsymbol{a}\rangle_t'=\big[\sum x_i a_i\big]_t'=\sum x_t'a_i=\langle \boldsymbol{x}', \boldsymbol{a}\rangle $, we have
$$
\langle \phi(\boldsymbol {b})-\phi(\boldsymbol {a}),\boldsymbol {b}-\boldsymbol {a}\rangle =\int_0^1 \langle v_t,\boldsymbol {b}-\boldsymbol {a}\rangle\,dt > (1-c)\int_0^1\|\boldsymbol {b}-\boldsymbol {a}\|^2\,dt = (1-c)\,\|\boldsymbol {b}-\boldsymbol {a}\|^2 > 0.
$$
In other words, the injectivity condition is satisfied:\ $\phi(\boldsymbol {a}) \ne \phi(\boldsymbol {b}) \ $for each point $ \boldsymbol {b} \in U \setminus \{ \boldsymbol {a} \}$. 
Lemma 3 is proved.   \hfill $\square$

\medskip
\textbf{Lemma 4.} {\it
The following inequality is true: $
|f^{(k+l)}_{x^ky^l}(x,y) - f^{(k+l)}_{x^ky^l}(0,0)| \le R(C_{k+1,l} + C_{k,l+1})\le 2MR.
$}

\medskip
{\it Proof.}
\begin{equation*}
\begin{split} 
|f^{(k+l)}_{x^ky^l}(x,y) - f^{(k+l)}_{x^ky^l}(0,0)| &= |f^{(k+l)}_{x^ky^l}(x,y) - f^{(k+l)}_{x^ky^l}(0,y) - f^{(k+l)}_{x^ky^l}(0,0)  + f^{(k+l)}_{x^ky^l}(0,y)|\le\\
&\le|f^{(k+l)}_{x^ky^l}(x,y) - f^{(k+l)}_{x^ky^l}(0,y)|+|f^{(k+l)}_{x^ky^l}(0,y) - f^{(k+l)}_{x^ky^l}(0,0)|=\\
&=\biggl| \int\limits_{0}^x f^{(k+l+1)}_{x^{k+1}y^l}(t,y)\,dt\biggr|+\biggl| \int\limits_{0}^y f^{(k+l+1)}_{x^ky^{l+1}}(0,t)\,dt\biggr|. 
\end{split} 
\end{equation*}
The lemma is proved. \hfill $\square$

\section
{Proof of Theorem 2} 
By Lemma 1, after the change $(x,y)\to(x_2,y_2)$, we have $f^{(a+b)}_{x_2^ay_2^b}({\bf 0})=0$ for all $a,b\in\mathbb Z_+$ with $a+b>0$, $a<3$ and $b<1+\frac k2$. By Lemma 2, $f=f(P)+x_2^3h(x_2,y_2) \pm y_2^{1+k/2}g(x_2,y_2)$ for some functions $h,g\in C^{\infty}(U_0)$. This immediately reduces the function to the required normal form $f=f(P)+{\tilde x}^3  \pm \tilde y^{1+k/2}$ by the coordinate change 
$\phi:\mathbb R_{x_2,y_2}^2\to\mathbb R_{\tilde x,\tilde y}^2 \  \mbox{with} \  
\tilde x = x_2 \sqrt[3]{h(x_2,y_2)}, \  
\tilde y  = y_2 ( g(x_2,y_2))^{\frac1{1+k/2}}$. From our proof of Lemma 2, we obtain explicite formulas for $h,g$ in terms of $S_{1+k/2}$ and $R_{1+k/2}$, namely:
$\displaystyle
f = \sum_{i=0}^{k/2}f_{y_2^i}^{(i)}(x_2,0)\frac{y_2^i}{i!} + R_{1+k/2}(x_2,y_2) 
= f(P) + x_2^3h(x_2,y_2) \pm y_2^{1+k/2}g(x_2,y_2)
= f(P) + {\tilde x}^3  \pm \tilde y^{1+k/2}$.

It remains to apply Lemma 3 to the coordinate transformation $ \phi $.
From the above formulas for the coordinate change $\phi$,
the bound $C_{\alpha\beta}<M$ and Lemma 4,
using and the Taylor expansion formula with a remainder in 
the Lagrange form or an integral remainder, we will obtain 
the required bound $ \|\phi'(\boldsymbol {x}) - I\| = \| A({\boldsymbol x})\| < C < 1$ for each point $\boldsymbol {x}\in U$. 
We will estimate the norm of the matrix in terms of its elements:
$
\|A\| \le \sqrt{\sum a_{ij}^2}{,}\quad  A = \{a_{ij}\}_{i,j=1}^n.
$

Let us proceed with detailed estimations, separately for the cases of $E_6$ and $E_8$.
By abusing notations, we will denote $(x_2,y_2)$ by $(x,y)$.

\subsection{The case of $E_6$} 
We compute the elements of the Jacobi matrix of 
$\phi$:
\begin{equation}\label{eq:i}
\frac{\partial \tilde x}{\partial x} = h^{\frac{1}{3}} + x\cdot \frac{1}{3} h^{-\frac{2}{3}}\cdot h_x',
\end{equation}
\begin{equation}\label{eq:ii}
\frac{\partial \tilde x}{\partial y} = x\cdot \frac{1}{3} h^{-\frac{2}{3}}\cdot h_y',
\end{equation}
\begin{equation}\label{eq:iii}
\frac{\partial \tilde y}{\partial x} = y\cdot \frac{1}{4} g^{-\frac{3}{4}} \cdot \biggl( \frac{1}{6} \int\limits_{0}^1 f_{xy^{4}}^{(5)}(x,sy)(1-s)^3\,ds \biggr),
\end{equation}
\begin{equation}\label{eq:iv}
\frac{\partial \tilde y}{\partial y} = g^{\frac{1}{4}} + y \cdot \biggl( \frac{1}{4}g^{-\frac{3}{4}} \cdot \frac{1}{6}\int\limits_{0}^1 s f_{y^{5}}^{(5)}(x,sy)(1-s)^3\,ds \biggr).
\end{equation}
If $M=0$, then $\phi=$\,id and everything is proved. Let further $M>0$, and therefore $R<1$. By Assumption 1, we have $h(0,0)=1$, $g(0,0)=1$. Thus, the Jacobi matrix at $\bf 0$ is the unit matrix $I$.

By using the above formulas (1)--(4) for Jacobi matrix' elements, 
let us estimate the elements of the Jacobi matrix $ \frac {\partial (\tilde x-x, \tilde y-y)} {\partial (x, y)} $  
and prove item (a) for the case of $E_6$. 

{\it Step 1.} Here we find an upper bound for $|1 - \frac{\partial \tilde x}{\partial x}|$. Remind that:
$$
\tilde x = x\sqrt[3]{h(x,y)} = \sqrt[3]{f(x,0) + yf_y'(x,0) + \frac{y^2}{2}f_{y^2}''(x,0) + \frac{y^3}{6}f_{y^3}'''(x,0)}.
$$
Denote $\tilde h = \tilde h(x,y) = h(x,y)-1$. 

Let us estimate $\tilde h(x,y)$. In the expression for $ h (x, y) $, we apply the Taylor expansion formula in $ x $ with a remainder in the Lagrange form to the coefficients of powers of $ y $:
\begin{equation*}
\begin{split}
|h(x,&y)-1| = \biggl|\frac{f_{x^3}'''(0,0)}{6} + x\frac{f_{x^4}^{(4)}(c_0,0)}{24} + \sum_{k=1}^3\big(\frac{f_{x^3}'''(c_k,0)}{6}\big)_{y^k}^{(k)}\frac{y^k}{k!}-1\biggr| = \\
&= \biggl|1 + x\frac{f_{x^4}^{(4)}(c_0,0)}{24} + y\frac{f_{x^3y}^{(4)}(c_1,0)}{6} + y^2\frac{f_{x^3y^2}^{(5)}(c_2,0)}{12} + y^3\frac{f_{x^3y^3}^{(6)}(c_3,0)}{36}-1\biggr| \le \\
&\le \frac{MR}{12}\biggl(\frac{5}{2}+R+\frac{R^2}{3}\biggr) < \frac{5}{24} \cdot \frac{MR}{1-R} \le \frac{5}{24} \cdot \frac{2M}{M+2}\cdot \frac{1}{1-2/(M+2)} = \frac{5}{24} \cdot 2 < \frac{1}{2}.
\end{split}
\end{equation*}
Hence $h(x,y) \in (0{.}5 , 1{.}5)$.
Therefore $\sqrt[3]{h(x,y)} \in (0{.}79, 1{.}21).$ 

Item (a) of Theorem 2 is proved for the function $h(x,y)$.
Set $c=0{.}5$. Then $|\tilde h(x,y)| < c$.

By the formula (\ref {eq:iii}) we have
$$ \frac{\partial \tilde x}{\partial x} = h^{\frac{1}{3}} + x\cdot \frac{1}{3} h^{-\frac{2}{3}}\cdot h_x' = (1 + \tilde h)^{\frac{1}{3}} + x\cdot \frac{1}{3} (1+ \tilde h)^{-\frac{2}{3}}\cdot h_x'.
$$
Let us estimate $|h_x'|$. We use Taylor' formula with an integral remainder:
\begin{equation*}
\begin{split} 
|h_x'| &= \biggl|\biggl(\sum_{k=0}^3 \frac{y^k}{k!}\frac{f_{y^k}^{(k)}(x,0)}{x^3}\biggr)_x'\biggr| =\biggl|\biggl(\sum_{k=0}^3 \frac{y^k}{k!}\frac{1}{2}\int\limits_{0}^1 f_{y^kx^3}^{(k+3)}(sx,0)(1-s)^2\,ds\biggr)_x'\biggr| = \\
&=\biggl|\sum_{k=0}^3 \frac{y^k}{k!\cdot2}\int\limits_{0}^1 sf_{y^kx^4}^{(k+4)}(sx,0)(1-s)^2\,ds\biggr| \le \sum_{k=0}^3 \frac{R^kM}{k!\cdot2}\int\limits_{0}^1 (s-2s^2+s^3)\,ds<\\
&< \frac{M}{24}\Big(1 + R + \frac{R^2}{2} + \frac{R^3}{6}\Big) < \frac1{24}\cdot\frac{M}{1-R} \le \frac1{12R}.
\end{split}
\end{equation*}
Set $c_x = \frac1{12R}$. 
We obtain a bound for $|\frac{\partial \tilde x}{\partial x} - 1|$ using the estimate of $\sqrt[3]{h}$:

\begin{equation*}
\begin{split} 
\Big|\frac{\partial \tilde x}{\partial x} - 1\Big| &= \big|\sqrt[3]{h} - 1 + x\cdot \tfrac{1}{3} (1+\tilde h)^{-\frac{2}{3}}h_x' \big|\le |\sqrt[3]{h} - 1| + \big|x\cdot \tfrac{1}{3} (1+\tilde h)^{-\frac{2}{3}}h_x'\big| < \\
& < 0{.}21 + \frac{c_xR}{3}\big(1-c)^{-\frac{2}{3}} = 0{.}21 + \frac{2^{\frac{2}{3}}}{3 \cdot 12}< 0{.}21 + 0{.}05=0{.}26.
\end{split}
\end{equation*} 

{\it Step 2.} Let us estimate $\frac{\partial \tilde x}{\partial y}$ from the formula (\ref {eq:ii}). \\
a) First we estimate $|h_y'|$ when $x\ne0$. We use Taylor' formula with an integral remainder:
\begin{equation*}
\begin{split}
|h_y'| &= \big|(f(x,0) + yf_y'(x,0) + \frac{y^2}{2}f_{y^2}''(x,0) + \frac{y^3}{6}f_{y^3}'''(x,0))_y'/x^3\big|=\\
&=\big|(f_y'(x,0) + yf_{y^2}''(x,0) + \frac{y^2}{2}f_{y^3}'''(x,0))/x^3\big|=\\
&=\biggl|\frac{1}{2}\int\limits_{0}^1 \biggl(f_{x^3y}^{(4)}(sx,0) + yf_{x^3y^2}^{(5)}(sx,0) +\frac{y^2}{2}f_{x^3y^3}^{(6)}(sx,0)\biggr)(1-s)^2\,ds\biggr| \le \\
&\le \frac{M}{6}(1+R+R^2) < \frac16\cdot\frac M{1-R} \le \frac1{3R}.
\end{split}
\end{equation*}
Set $c_y = \frac1{3R}$.

b) In the formula for $|\frac{\partial \tilde x}{\partial y}|$ we expand $ h ^ {- \frac{2}{3}} $ in a Taylor series in $ \tilde h $ and take $ c_{\tilde h} = c_{\tilde h}(x, y) \in [0,\tilde h (x, y)] $ from the Taylor-Lagrange formula:

\begin{equation*}
\begin{split} 
\big|\frac{\partial \tilde x}{\partial y}\big| &= \big|x\cdot \tfrac{1}{3} \big(1+\tilde h\big)^{-\frac{2}{3}}\cdot h_y'\big|= \big|\frac{x}{3}\big(1-\tfrac{2}{3}\tilde h + \tfrac{5}{9}(1+c_{\tilde h})^{-\frac{8}{3}}{\tilde h}^2\big)h_y'\big|\le \\
&\le \frac{R}{3}\big(1 + \tfrac{2}{3}c + \tfrac{5}{9}(1+c_{\tilde h})^{-\frac{8}{3}}c^2\big)c_y < \frac{c_yR}{3}\big(1 + \tfrac{2}{3}c + \tfrac{5}{9}(1-c)^{-\frac{8}{3}}c^2\big)=\\
&=\frac 1{9}\big(1 + \tfrac1{3} + \tfrac{5}{9}2^{\frac{8}{3}}/4\big) < 0{.}25.
\end{split}
\end{equation*}

{\it Step 3.} Let us estimate $ \frac{\partial \tilde y}{\partial x} $ from the formula (\ref {eq:iii}). Let us first estimate separately the factors of this expression.

a) Auxiliary Assessment. By Lemma 4, we have	
\begin{equation*}
\begin{split} 
\biggl|\int\limits_{0}^1 \biggl( f_{y^4}^{(4)}(x,sy) - f_{y^4}^{(4)}(0,0)\biggr)(1-s)^3\,ds\biggr| &\le \int\limits_{0}^1 2MR(1-s)^3\,ds = \frac{MR}{2} <1.
\end{split}
\end{equation*}

b) By Assumption 1, we have $f_{y^4}^{(4)}(0,0) = 24$. We get the following lower bound:
\begin{equation*}
\begin{split} 
\biggl|\int\limits_{0}^1 f_{y^4}^{(4)}(x,sy)(1-s)^3\,ds\biggr| &= \biggl|\int\limits_{0}^1 \biggl( f_{y^4}^{(4)}(x,sy) - f_{y^4}^{(4)}(0,0) + f_{y^4}^{(4)}(0,0)\biggr)(1-s)^3\,ds\biggr|\\
& \geq 6 - \frac{MR}{2} > 5 > 0 .
\end{split}
\end{equation*}
Item (a) of Theorem 2
has been completely proved for $E_6$. 
Set $\tilde c = 1 - \frac{MR}{12}$. 

Then 
\begin{equation*}
\begin{split} 
\big|\frac{\partial \tilde y}{\partial x}\big| &= \biggl|y \frac{1}{4} \biggl( \frac{1}{6} \int\limits_{0}^1 f_{y^4}^{(4)}(x,sy)(1-s)^3\,ds \biggr) ^{-\frac{3}{4}} \biggl( \frac{1}{6} \int\limits_{0}^1 f_{xy^4}^{(5)}(x,sy)(1-s)^3\,ds \biggr)\biggr|\le\\
&\le R\frac{1}{4}\tilde c^{-\frac{3}{4}}\frac{M}{24}<\frac{MR}{96}\biggl(1-\frac{MR}{12}\biggr)^{-\frac{3}{4}}<\frac{1}{48}\biggl(1-\frac{1}{6}\biggr)^{-\frac{3}{4}}<0{.}03.
\end{split}
\end{equation*}

{\it Step 4.} Let us estimate $|1 - \frac{\partial \tilde y}{\partial y}|$ from formula (\ref {eq:iv}). 
\begin{equation*}
\begin{split}
\Big|\frac{\partial \tilde y}{\partial y} - 1 \Big|&=\Biggl|\sqrt[4]{\frac{1}{6} \int_{0}^1 f_{y^4}^{(4)}(x,sy)(1-s)^3\,ds } - 1 + \\
&+\frac{y}{4}\biggl( \frac{1}{6}\int\limits_{0}^1 f_{y^4}^{(4)}(x,sy)(1-s)^3\,ds \biggr) ^{-\frac{3}{4}}\frac{1}{6}\int\limits_{0}^1 s f_{y^5}^{(5)}(x,sy)(1-s)^3\,ds\Biggr|\le\\
&<\Big(1 + \frac{MR}{12}\Big)^{\frac{1}{4}} - 1 + \frac{MR}{80}\Big(1-\frac{MR}{12}\Big)^{-\frac{3}{4}}<\\
&<\Big(1 + \tfrac{1}{6}\Big)^{\frac{1}{4}} - 1 + \tfrac{1}{40}\Big(1-\tfrac{1}{6}\Big)^{-\frac{3}{4}}<0{.}07.
\end{split}
\end{equation*}

For each point $\boldsymbol{x}\in U$, we have: \\ $ \|\phi'(\boldsymbol {x}) - I\| = \| A({\boldsymbol x})\| \le \sqrt{0{.}26^2+0{.}25^2+0{.}03^2+0{.}07^2}= \sqrt{0{.}1359} < 0{.}4 = C < 1$, that proves item (b) of the Theorem 2 for $C=\frac25$. By Lemma 3 the coordinate change $\phi$ is injective in $U$, that proves item (c), except for the properties of $ \phi (U) $.

From the last assertion of Lemma 3 and [4, Cor.~8.3, Step 1], we conclude that $\phi(U)$ contains the open disk of radius $ (1 - C) R $ centred at the origin. This completes our proof of Theorem 2 for the case of a singularity $E_6$.

\subsection{The case of $E_8$}
Consider the coordinate change  $ \tilde x = x\sqrt[3]{h(x,y)} = \sqrt[3]{\sum_{k=0}^{4} \frac{y^k}{k!}f_{y^k}^{(k)}(x,0)}$, $ \tilde y = y\sqrt[5]{g(x,y)}= y \left[ \frac{1}{24}\int_{0}^1 f_{y^{5}}^{(5)}(x,sy)(1-s)^4\,ds\right]^{\frac{1}{5}}$.

{\it Step 1.} As in the case of $E_6$, we find a bound for the term $\tilde h(x,y)$, in order to estimate $|1 - \frac{\partial \tilde x}{\partial x}|$. The difference from the case of $E_6$ is an additional term in the sum of absolute values that can be bound by $\frac{MR^4}{144}$, then the sum can be bound by the same geometric progression. In detail:
\begin{equation*}
\begin{split}
|h(x,y)-1| &= \biggl|\frac{f_{x^3}'''(0,0)}{6} + x\frac{f_{x^4}^{(4)}(c_0,0)}{24} + \sum_{k=1}^4\big(\frac{f_{x^3}'''(c_k,0)}{6}\big)_{y^k}^{(k)}\frac{y^k}{k!}-1\biggl| \le \\
&\le \frac{RC_{40}}{24} + \frac{RC_{31}}{6} + \frac{R^2C_{32}}{12} + \frac{R^3C_{33}}{36} + \frac{R^4C_{34}}{144} \le \frac{MR}{12}\biggl(\frac{5}{2}+R+\frac{R^2}{3}+\frac{R^3}{12}\biggr) < \frac{5}{24} \cdot 2< \frac{1}{2}.
\end{split}
\end{equation*}
The bound from before holds: $h(x,y) \in (0{.}5 , 1{.}5)$ or $|\tilde h| < 0{.}5$. Then $\sqrt[3]{h(x,y)} \in (0{.}79, 1{.}21)$.

Similar arguments work for the estimation of the term $|h_x'|$. Thus we can bound it as we did it in the case of $E_6$:
$$
|h_x'| <\frac1{12R}.
$$
Set $c_x = \frac{1}{12R}$.

The element $\frac{\partial \tilde x}{\partial x}$ has the same representation in terms of the function $h$ as in the case $E_6$ and the same estimation for every term. Hence the following bound holds:
$$
\Big|\frac{\partial \tilde x}{\partial x} - 1\Big| < 0{.}26.
$$

{\it Step 2.} a) For estimating the term $|h_y'|$ when $x\ne0$, we can use the bound from the case of $E_6$:
\begin{equation*}
|h_y'| = \big|\biggl(\sum_{k=0}^4 \frac{y^k}{k!}\frac{f_{y^k}^{(k)}(x,0)}{x^3}\biggr)_y'\big| \le \frac{1}{3R}.
\end{equation*}

b) The element $\frac{\partial \tilde x}{\partial y}$ has the same representation in terms of the function $h$ as in the case of $E_6$. Hence the following bound holds:
$$
\Big|\frac{\partial \tilde x}{\partial y} \Big| < 0{.}25.
$$

{\it Step 3.} In this step, we have orders of partial derivatives and a constant $\tilde c$ that differ from the case of $E_6$.

a) Auxiliary Assessment. By Lemma 4, we have 
\begin{equation*}
\begin{split} 
\bigg|\int_{0}^1 \biggl( f_{y^5}^{(5)}(x,sy) - f_{y^5}^{(5)}(0,0)\biggr)(1-s)^4\,ds\bigg| &\le \int_{0}^1 R(C_{15}+C_{06})(1-s)^4\,ds = \frac{(C_{15}+C_{06})R}{5} \le \frac{2MR}{5}.
\end{split}
\end{equation*}

b) By Assumption 1, we have $f_{y^5}^{(5)}(0,0) = 120$. We get the following lower bound:
\begin{equation*}
\begin{split} 
\biggl|\int_{0}^1 f_{y^5}^{(5)}(x,sy)(1-s)^4\,ds\biggl| &=  \biggl|\int_{0}^1 \biggl( f_{y^5}^{(5)}(x,sy) - f_{y^5}^{(5)}(0,0) + f_{y^5}^{(5)}(0,0)\biggr)(1-s)^4\,ds\biggl| \geq 24 - \frac{2MR}{5}.
\end{split}
\end{equation*}
Set $\tilde c = 1 - \frac{MR}{60}$. Then
\begin{equation*}
\begin{split} 
\big|\frac{\partial \tilde y}{\partial x}\big| &= \biggl|y\cdot \frac{1}{5} \biggl( \frac{1}{24} \int_{0}^1 f_{y^5}^{(5)}(x,sy)(1-s)^4\,ds \biggr) ^{-\frac{4}{5}} \cdot \biggl( \frac{1}{24} \int_{0}^1 f_{xy^5}^{(6)}(x,sy)(1-s)^4\,ds \biggr)\biggl|\le\\
&\le R\frac{1}{5}\tilde c^{-\frac{4}{5}}\frac{C_{15}}{120}<\frac{MR}{600}\biggl(1-\frac{MR}{60}\biggr)^{-\frac{4}{5}}<\frac{1}{300}\biggl(1-\frac{1}{30}\biggr)^{-\frac{4}{5}}<0{.}004<0{.}03.
\end{split}
\end{equation*}

{\it Step 4.} Let us estimate $|1 - \frac{\partial \tilde y}{\partial y}|$:
\begin{equation*}
\begin{split}
\Big|\frac{\partial \tilde y}{\partial y} - 1 \Big|&=\Bigg|\sqrt[5]{\frac{1}{24} \int_{0}^1 f_{y^5}^{(5)}(x,sy)(1-s)^4\,ds } - 1 + \\
&+ \frac{y}{5}\biggl( \frac{1}{24}\int_{0}^1 f_{y^5}^{(5)}(x,sy)(1-s)^4\,ds \biggr) ^{-\frac{4}{5}}\frac{1}{24}\int_{0}^1 s f_{y^6}^{(6)}(x,sy)(1-s)^4\,ds\Bigg|\le\\
&\le \bigg|\bigg(\frac{1}{24}\int_{0}^1 \biggl( f_{y^5}^{(5)}(x,sy) - f_{y^5}^{(5)}(0,0) + f_{y^5}^{(5)}(0,0)\biggr)(1-s)^4\,ds\bigg)^{\frac{1}{5}} - 1\bigg| + |\dots|\le \\
& \le \Big|\Big(1 - \frac{(C_{15}+C_{06})R}{120}\Big)^{\frac{1}{5}} - 1\Big| + R\tilde c^{-\frac{4}{5}}\frac{C_{06}}{5\cdot 24 \cdot 30}<1-\Big(1 - \frac{MR}{60}\Big)^{\frac{1}{5}} + \frac{MR}{3600}\Big(1-\frac{MR}{60}\Big)^{-\frac{4}{5}}<\\
&<1-\Big(1 - \tfrac{1}{30}\Big)^{\frac{1}{5}} + \tfrac{1}{1800}\Big(1-\tfrac{1}{30}\Big)^{-\frac{4}{5}}< 0{.}008<0{.}07.
\end{split}
\end{equation*}

Thus we found the required bound for a singularity of the type $E_8$: $\|\phi'(\boldsymbol {x}) - I\|<\sqrt{0{.}26^2+0{.}25^2+0{.}004^2+0{.}008^2}<0{.}4<1$.

Theorem 2 is completely proved. \hfill $\square$

\section{The case of $A_n$} 

\medskip
\textbf{Definition.} {\it
A smooth function $f(u_1, u_2) $ has a singularity $A_n$ {\rm($n\ge1$)} at its critical point $P\in\mathbb R^2$ if
\begin{itemize}
\item[\rm(i)] the first differential $df(P) = 0$ and the second differential $d^2f(P) \ne 0$ and is a perfect square;
\item[\rm(ii)] some codition on coefficients of the Taylor series of $f$ at $P$ holds (see below).
\end{itemize}
Consider a linear change $(u_1,u_2)\to(x,y)$ such that $d^2f(P)=\pm(dx)^2$. By using a non-linear change $(x,y)\to(x_1=x+y^2Q(y),y_1=y)$ for some polynomial $Q(y)$, one can achive that $ f^{(a+b)}_{x_1^ay_1^b}(P)=0$ for all $a,b\in\mathbb Z_+$ with $0<\frac{a}{2}+\frac{b}{n+1}<1$. 
Notice that $f''_{x_1^2}(P) = \pm 2$. The condition on the Taylor coefficients is as follows: $f^{(n+1)}_{y_1^{n+1}}(P) \ne 0$.}

We will assume that $P=(0,0)={\bf 0}$ in the coordinates $x, y$.

\medskip
\textbf{Assumption 2.} For singularities $A_n$ {\rm($n\ge1$)}, assume that
$f_{y_1^{n+1}}^{(n+1)}(0,0) = \pm(n+1)!$.

\medskip
\textbf{Theorem 3} (Estimating the radius of a neighborhood for the coordinate change)\textbf{.} {\it 
In the case of a singularity $A_n$ ($n\ge1$), under the hypotheses of Assumption 2, suppose that, in a neighbourhood 
$U_0 =\{(x_1,y_1)\mid \max(|x_1|,|y_1|) < R_0\}$ of $\bf 0$, the following estimates hold:
$C_{\alpha\beta}=\sup_{U_0}\big|f_{x_1^{\alpha}y_1^{\beta}}^{(\alpha+\beta)}(x_1,y_1)\big|\le M$ for all $(\alpha,\beta)\in\{ (0,n+2), (0,n+1), (2,i), (3,0), (3,i) \mid i = 0, 1, \dots, k\}$, where $R_0>0$, $M\ge0$. 
Then, in the neighborhood $ U = \{(x_1, y_1) \mid \max (| x_1 |, | y_1 |) <R \} $, with $R =\min\{R_0,\frac{3}{4M+3}\}$, there is a coordinate change of the form
$\phi:(x_1,y_1)\to(\tilde x = x_1\sqrt{h(x_1,y_1)}, \tilde y=y_1(g(x_1,y_1))^{\frac1{k+1}})$ that reduces $ f $ to the normal form $ f = f(P) \pm \tilde x^2 \pm \tilde y^{n+1}$ of $A_n$. In detail: the coordinate change $\phi|_U$ satisfies the conditions {\rm (a), (b), (c)} from Theorem 2 with $C=0{.}93$.}

\medskip
{\it Proof.} By abusing notations, we will denote $(x_1,y_1)$ by $(x,y)$.
Consider the coordinate chage  
$$ 
\tilde x = x\sqrt{h(x,y)} = \sqrt{\sum_{k=0}^{n} \frac{y^k}{k!}f_{y^k}^{(k)}(x,0)},\quad \tilde y = y\sqrt[n+1]{g(x,y)}= y \left[ \frac{1}{n!}\int_{0}^1 f_{y^{n+1}}^{(n+1)}(x,sy)(1-s)^{n}\,ds\right]^{\frac{1}{n+1}}.
$$
We extend our proof of Theorem 2 (about singularities $E_6, E_8$) to the case of singularities $A_n$.

\medskip
{\it Step 1.} 
a) Let us estimate $\tilde h(x,y)$:
\begin{equation*}
\begin{split}
|h(x,y)-1| &= \biggl|\frac{f_{x^2}''(0,0)}{2} + x\frac{f_{x^3}'''(c_0,0)}{6} + \sum_{k=1}^{n}\frac{f_{x^2y^k}^{(k+2)}(c_k,0)}{2}\frac{y^k}{k!}-1\biggl| \le \\
&\le \frac{RC_{30}}{6}+\sum_{k=1}^{n}\frac{C_{2k}R^k}{2\cdot k!} \le \frac{2MR}{3(1-R)}\le\frac{1}{2}.
\end{split}
\end{equation*}
Thus $h(x,y) \in (0{.}5 , 1{.}5)$, or $|\tilde h| < 0{.}5 =: c$. Then $\sqrt{h(x,y)} \in (0{.}7, 1{.}3)$. 

We have
$$ 
\frac{\partial \tilde x}{\partial x} = h^{\frac{1}{2}} + x\cdot \frac{1}{2} h^{-\frac{1}{2}}\cdot h_x' = (1 + \tilde h)^{\frac{1}{2}} + x\cdot \frac{1}{2} (1+ \tilde h)^{-\frac{1}{2}}\cdot h_x'.
$$
Let us estimate $|h_x'|$:
\begin{equation*}
\begin{split} 
|h_x'| &= \biggl|\biggl(\sum_{k=0}^{n} \frac{y^k}{k!}\frac{f_{y^k}^{(k)}(x,0)}{x^2}\biggr)_x'\biggl| =\biggl|\biggl(\sum_{k=0}^{n} \frac{y^k}{k!}\int_{0}^1 f_{y^kx^2}^{(k+2)}(sx,0)(1-s)\,ds\biggr)_x'\biggl| = \\
&=\biggl|\sum_{k=0}^{n} \frac{y^k}{k!}\int_{0}^1 sf_{y^kx^3}^{(k+3)}(sx,0)(1-s)\,ds\biggl| \le \sum_{k=0}^{n} \frac{R^kC_{3k}}{k!}\int_{0}^1 (s-s^2)\,ds\le\\
&\le \frac{M}{6(1-R)}\le \frac{1}{8R}.
\end{split}
\end{equation*}
Set $c_x = \frac{1}{8R}$.
Let us estimate $|\frac{\partial \tilde x}{\partial x} - 1|$:
\begin{equation*}
\begin{split} 
\Big|\frac{\partial \tilde x}{\partial x} - 1\Big| &= \big|\sqrt{h} - 1 + x\cdot \tfrac{1}{2} (1+\tilde h)^{-\frac{1}{2}}h_x' \big|\le |\sqrt{h} - 1| + \big|x\cdot \tfrac{1}{2} (1+\tilde h)^{-\frac{1}{2}}h_x'\big| < \\
& < 0{.}3 + \frac{c_xR}{2}\big(1-c)^{-\frac{1}{2}} < 0{.}3 + \frac{2^{\frac{1}{2}}}{2 \cdot 8}< 0{.}3 + 0{.}1=0{.}4.
\end{split}
\end{equation*} 

{\it Step 2.} a) Let us estimate $|h_y'|$ when $x\ne0$:
\begin{equation*}
\begin{split}
|h_y'| &= \big|\biggl(\sum_{k=0}^n \frac{y^k}{k!}\frac{f_{y^k}^{(k)}(x,0)}{x^2}\biggr)_y'\big|=\big|\biggl(\sum_{k=0}^{n-1} \frac{y^k}{k!}\frac{f_{y^{k+1}}^{(k+1)}(x,0)}{x^2}\biggr)\big|=\\
&=\biggl|\int_{0}^1 \biggl(\sum_{k=0}^{n-1}\frac{y^k}{k!}f_{x^2y^{k+1}}^{(k+3)}(sx,0)\biggr)(1-s)\,ds\biggl| \le \frac{1}{2}\sum_{k=0}^{n-1}\frac{R^kC_{2,k+1}}{k!} \le \frac{M}{2(1-R)}\le \frac{3}{8R}.
\end{split}
\end{equation*}
Set $c_y = \frac{3}{8R}$.\\
b) Let us estimate $|\frac{\partial \tilde x}{\partial y}|$:
\begin{equation*}
\begin{split} 
\big|\frac{\partial \tilde x}{\partial y}\big| &= \big|x\cdot \tfrac{1}{2} \big(1+\tilde h\big)^{-\frac{1}{2}}\cdot h_y'\big|= \big|\frac{x}{2}\big(1-\tfrac{1}{2}\tilde h + \tfrac{3}{8}(1+c_{\tilde h})^{-\frac{5}{2}}{\tilde h}^2\big)h_y'\big|\le \\
&\le \frac{R}{2}\big(1 + \tfrac{1}{2}c + \tfrac{3}{8}(1+c_{\tilde h})^{-\frac{5}{2}}c^2\big)c_y < \frac{c_yR}{2}\big(1 + \tfrac{1}{2}c + \tfrac{3}{8}(1-c)^{-\frac{5}{2}}c^2\big)=\\
&=\frac {3}{16}\big(1 + \tfrac1{4} + \tfrac{3}{8}2^{\frac{5}{2}}/4\big) < 0{.}4.
\end{split}
\end{equation*}

{\it Step 3.} a) Auxiliary Assessment. By Lemma 4:	
\begin{equation*}
\begin{split} 
\bigg|\int_{0}^1 \biggl( f_{y^{n+1}}^{(n+1)}(x,sy) - f_{y^{n+1}}^{(n+1)}(0,0)\biggr)(1-s)^n\,ds\bigg| &\le \int_{0}^1 R(C_{1,n+1}+C_{0,n+2})(1-s)^n\,ds \le \frac{2MR}{n+1}.
\end{split}
\end{equation*}
b) By assumption 3 we have $f_{y^{n+1}}^{(n+1)}(0,0) = (n+1)!$. We get the following lower bound:
\begin{equation*}
\begin{split} 
\biggl|\int_{0}^1 f_{y^{n+1}}^{(n+1)}(x,sy)(1-s)^n\,ds\biggl| &= \biggl|\int_{0}^1 \biggl( f_{y^{n+1}}^{(n+1)}(x,sy) - f_{y^{n+1}}^{(n+1)}(0,0) + f_{y^{n+1}}^{(n+1)}(0,0)\biggr)(1-s)^n\,ds\biggl| \geq n! - \frac{2MR}{n+1}.
\end{split}
\end{equation*}
Set $\tilde c = 1 - \frac{2MR}{(n+1)!}$. Then the following sequence of inequalities holds:
\begin{equation*}
\begin{split} 
\big|\frac{\partial \tilde y}{\partial x}\big| &= \biggl|y\cdot \frac{1}{n+1} \biggl( \frac{1}{n!} \int_{0}^1 f_{y^{n+1}}^{(n+1)}(x,sy)(1-s)^n\,ds \biggr) ^{-\frac{n}{n+1}} \cdot \biggl( \frac{1}{n!} \int_{0}^1 f_{xy^{n+1}}^{(n+2)}(x,sy)(1-s)^n\,ds \biggr)\biggl|\le\\
&\le R\frac{1}{n+1}\tilde c^{-\frac{n}{n+1}}\frac{C_{1,n+1}}{(n+1)!}\le\frac{MR}{(n+1)(n+1)!}\biggl(1-\frac{2MR}{(n+1)!}\biggr)^{-\frac{n}{n+1}}<\frac{3}{4\cdot 4}\biggl(1-\frac{3}{2(n+1)!}\biggr)^{-\frac{n}{n+1}}\le \frac38.
\end{split}
\end{equation*}
Let us prove the latter inequality in this sequence. For $n=1$ this inequality is in fact an equality. For $n>1$ we have $\tilde c^{-\frac{n}{n+1}}<\tilde c^{-1}$, moreover the function $\tilde c^{-1}|_{MR=\frac3{4}}$ of $n$ is monotone decreasing. Since $\tilde c^{-1}|_{MR=\frac3{4}} < 2$ for $n = 2$, we have the same for all $n>2$. 
This proves the desired inequality.

{\it Step 4.} Let us estimate $|1 - \frac{\partial \tilde y}{\partial y}|$:
\begin{equation*}
\begin{split}
\Big|\frac{\partial \tilde y}{\partial y} - 1 \Big|&=\Bigg|\sqrt[n+1]{\frac{1}{n!} \int_{0}^1 f_{y^{n+1}}^{(n+1)}(x,sy)(1-s)^n\,ds } - 1 + \\
&+\frac{y}{n+1}\biggl( \frac{1}{n!}\int_{0}^1 f_{y^{n+1}}^{(n+1)}(x,sy)(1-s)^n\,ds \biggr) ^{-\frac{n}{n+1}}\frac{1}{n!}\int_{0}^1 s f_{y^{n+2}}^{(n+2)}(x,sy)(1-s)^n\,ds\Bigg|\le\\
&\le \bigg|\bigg(\frac{1}{n!}\int_{0}^1 \biggl( f_{y^{n+1}}^{(n+1)}(x,sy) - f_{y^{n+1}}^{(n+1)}(0,0) + f_{y^{n+1}}^{(n+1)}(0,0)\biggr)(1-s)^n\,ds\bigg)^{\frac{1}{n+1}} - 1\bigg| + |\dots|\le \\
& \le \Big|\Big(1 - \frac{(C_{1,n+1}+C_{0,n+2})R}{(n+1)!}\Big)^{\frac{1}{n+1}} - 1\Big| + R\tilde c^{-\frac{n}{n+1}}\frac{C_{0,n+2}}{(n+2)(n+1)(n+1)!}<\\
&<1-\Big(1 - \frac{2MR}{(n+1)!}\Big)^{\frac{1}{n+1}} + \dots = 1-\tilde c^\frac1{n+1} + \frac{MR\tilde c^{-\frac{n}{n+1}}}{(n+2)(n+1)(n+1)!}  <\\
& < \frac3{2\cdot (n+1)!} + \frac{3\cdot 2}{4\cdot 3\cdot 2\cdot 2}\le\frac1{2}\le \frac58
\end{split}
\end{equation*}
if $n>1$. And when $n=1$, we get $1-\tilde c^\frac12 \le \frac12$, therefore $|\frac{\partial \tilde y}{\partial y} - 1|<\frac12+\frac18=\frac58$.

Thus for $A_n$ we have $\|\phi'(\boldsymbol {x}) - I\|<\sqrt{0{.}4^2+0{.}4^2+\frac3{8}^2+\frac58^2}<0{.}93<1$ for all $n\ge1$.

Theorem 3 is proved. \hfill $\square$


\end{otherlanguage}



\end{document}